\documentclass[matann, numbook, envcountsame, final]{svjour}
\usepackage{latexsym, amssymb, theorem, amsmath, amscd}

\def\CM{Cohen-Macaulay }
\def\tor{\mathrm{Tor}}
\def\rank{\mathop{\mathrm{rank}}}
\def\dim{\mathop{\mathrm{dim}}}
\def\ch{\mathrm{ch}}
\def\Supp{\mathop{\mathrm{Supp}}}
\def\Spec{\mathop{\mathrm{Spec}}}
\def\div{\mathrm{\mathop{div}}}
\def\chinf{\mathop{\chi_{\infty}}}
\def\Fdot{F_{\bullet}}
\def\Gdot{G_{\bullet}}

\begin{document}

\title{Intersection multiplicities over Gorenstein rings}

\author{Claudia M. Miller \and Anurag K. Singh\thanks{\emph{Present
address:} Department of Mathematics, University of Utah, 155 South, 1400 East,
Salt Lake City, UT 84112--0090, USA}}   
                  
\institute{Department of Mathematics, The University of Michigan, East Hall, \\
525 East University Avenue, Ann Arbor, MI 48109--1109, USA \\
E-mail: {\tt clamille@math.lsa.umich.edu,  singh@math.utah.edu}}

\dedication{Dedicated to Professor Jack E. McLaughlin}

\date{Received: February 11, 1999}

\maketitle

\begin{abstract}
Let $R$ be a complete local ring of dimension $d$ over a perfect field of prime
characteristic $p$, and let $M$ be an $R$-module of finite length and
finite projective dimension. S.~Dutta showed that the equality
$\lim_{n\to\infty}\frac{\ell(F^n_R(M))}{p^{nd}} = \ell(M)$ holds when the ring 
$R$ is a complete intersection or a Gorenstein ring of dimension at most $3$.
We construct a module over a Gorenstein ring $R$ of dimension five for which
this equality fails to hold. This then provides an example of a nonzero Todd 
class $\tau_3(R)$, and of a bounded free complex whose local Chern character 
does not vanish on this class.   
\end{abstract}

\section{Introduction}
\label{intro}

Let $R$ be a complete local ring of dimension $d$ over a perfect field of
prime characteristic $p$, and let $\Gdot$ be a bounded complex of free modules
with finite length homology. In \cite{Frobmult} S.~Dutta introduced the 
{\it limit multiplicity}
$$ 
\chinf(\Gdot) = \lim_{n \to \infty}\frac{\chi(F^n(\Gdot))}{p^{nd}}, 
$$ 
where $F^n(-)$ denotes the $n$\,th iteration of the Frobenius functor. 
C.~Szpiro and L.~Peskine showed in \cite{PS2} that the equality
$$ 
\chi(F^n(\Gdot)) = p^{nd} \chi(\Gdot)
$$ 
holds for a graded complex over a graded ring, and Szpiro conjectured that
this would hold in general over a \CM ring, see \cite{Szpiro}. In the specific
case that $\Gdot$ is the resolution of a module $M$ of finite length and finite
projective dimension, the conjecture then asserts that 
$$ 
\ell(F^n(M)) = p^{nd}\ell(M).  
$$ 
Dutta showed that this equality does hold if the ring $R$ is a complete 
intersection or a Gorenstein ring of dimension at most $3$, see \cite{Frobmult}. 
In \cite{RobMSRI} 
P.~Roberts constructed a counterexample to Szpiro's conjecture over a \CM ring
of dimension three using the famous example of negative Serre intersection
multiplicity due to Dutta, Hochster and McLaughlin, \cite{DHM}. The question
however remained open for Gorenstein rings, and the main aim of this paper is
to demonstrate that the conjecture is false in general over Gorenstein 
rings. Specifically, we construct a module $M$ of finite length and finite 
projective dimension over a Gorenstein ring $R$ of dimension five such that
$$
\lim_{n\to\infty}\frac{\ell(F^n_R(M))}{p^{5n}} \neq \ell(M).
$$

Using techniques similar to those used in \cite{DHM}, we first construct a 
module $N$ of finite length and of finite projective dimension over the 
hypersurface
$$ 
A=K[U,V,W,X,Y,Z]_m/(UX+VY+WZ), 
$$ 
where $m$ is the ideal $(U,V,W,X,Y,Z)$, such that $\chi(N, A/P) = -2$,
where $P$ is the prime ideal $(u,v,w)$. (Lower case letters will
denote the images of the corresponding variables.) We believe this is
of independent interest since it is an example of two modules, one of
whom has finite projective dimension, with a nonvanishing intersection
multiplicity where the sum of the dimensions of the modules is $\dim A
- 2$. It should be pointed out that if $N_1$ and $N_2$ are two
modules, each of which has finite projective dimension over $A$, then
the condition
$$
\dim N_1 + \dim N_2 < \dim A
$$ 
does imply $\chi(N_1, N_2) = 0$ by the main theorems of \cite{Robmult} or 
\cite{GS}.

The limit multiplicity has an interpretation in terms of localized
Chern characters and the local Riemann-Roch formula, see
\cite{RobMSRI} or \cite{RobBook}. The results mentioned above then
provide an example of a nonzero Todd class $\tau_3(R)$ over a
Gorenstein ring $R$ of dimension five, and of a bounded free complex
whose local Chern character does not vanish on this class. In
\cite{Kurano} K.~Kurano does obtain a Gorenstein ring $R$ of dimension
five with a nonzero Todd class $\tau_3 (R)$, but it is not known if
there is a free complex whose local Chern character does not
vanish on this Todd class.

\section{Background}
\label{background}

In this section we give a brief summary of the relevant terminology as well as
record some well-known facts about $\chi$ and $\chi_{\infty}$ that we shall 
find useful later in our work.

For two modules $M$ and $N$ over a local ring $(R,m)$ such that $M\otimes_R
N$ has finite length and $M$ has finite projective dimension, the Serre
intersection multiplicity is defined as 
$$
\chi(M,N) = \sum_{i=0}^{\dim R}(-1)^i \ell(\tor_i^R(M,N))
$$ 
where $\ell(-)$ denotes the length. This definition does agree with the 
geometric notion of intersection multiplicity, see \cite{Serre}.
For a bounded complex $\Gdot$
$$
0 \to G_s  \to \cdots \to G_1 \to G_0 \to 0
$$ 
with finite length homology, we define
$$
\chi(\Gdot) = \sum_{i=0}^s (-1)^i \ell(H_i(\Gdot)).
$$ 

Let $R$ be a ring of positive characteristic $p$ and
dimension $d$. Using the Frobenius endomorphism $f$ of $R$ (which
takes $r$ to $r^p$ for $r\in R$), Peskine and Szpiro defined the
Frobenius functor $F_R(-)$. This functor takes an $R$-module $M$ to
$$
F_R(M)=M\otimes_R{}^f R,
$$ 
where ${ }^f R$ is $R$ viewed as a module over itself with a left action via
the Frobenius endomorphism, and a right action via the usual multiplication. 
We use $F_R^n(-)$ (or simply $F^n(-)$) to denote the $n$\,th iteration of the
Frobenius functor. For a bounded free complex $\Gdot$ with finite length 
homology, Dutta defined
$$
\chinf(\Gdot) = \lim_{n\to\infty}
\frac{\chi(F^n(\Gdot))}{p^{nd}}.
$$
If $\Gdot$ is a finite free resolution of an $R$-module $M$ of finite length 
and finite projective dimension, then
$$
\chi(\Gdot) = \ell(M) {\textrm{\qquad and\qquad}}
\chinf(\Gdot) = \lim_{n\to\infty}\frac{\ell(F^n_R(M))}{p^{nd}}.
$$
For the second assertion note that 
$$
H_i(F^n_R(\Gdot)) = \tor_i^R(M,{ }^{f^n} R) = 0
$$ 
for all $i\ge 1$ and $n\geq 0$, by \cite[Theorem 1.7]{PS1}.

\begin{lemma}
\label{rank}
Let $R$ be an integral domain of characteristic $p>0$ and let $S$ be a 
module-finite extension ring which has rank $r$ as an $R$-module. Then for any 
bounded complex $\Gdot$ of free $R$-modules which has finite length homology, 
$$
\chinf(\Gdot\otimes_R S) = r\cdot\chinf(\Gdot).
$$
\end{lemma}

\begin{proof}
Since $S$ is an $R$-module of rank $r$, there exists a short exact sequence
$$
\begin{CD}
0 @>>> R^r @>>> S @>>> T @>>> 0
\end{CD}
$$
where $T$ is a torsion $R$-module and thus has dimension less than $d$. Since 
$\chi(F^n_R(\Gdot)\otimes_R -)$ is additive on short exact sequences, we get
$$
\chi(F^n_R(\Gdot)\otimes_R S) 
= r\cdot\chi(F^n_R(\Gdot)) + \chi(F^n_R(\Gdot)\otimes_R T). {\textrm \qquad (*)}
$$
Since the module $T$ has dimension less than the dimension of $R$, 
$$
\lim_{n\to\infty}\frac{\chi(F^n_R(\Gdot)\otimes_R T)}{p^{nd}} = 0
$$
by \cite[Proposition 1]{Seibert}. We obtain the desired equality by dividing 
the equation $(*)$ by $p^{nd}$ and forming the appropriate limits since
$$
F^n_R(\Gdot)\otimes_R S = F^n_S(\Gdot\otimes_R S). 
$$
\qed
\end{proof}

We will also use the following lemma, which is patterned after
Proposition 2.4 of \cite{DHM}.

\begin{lemma}
\label{projdim}
Let $(R,m,K)$ be a local ring, and let $M$ be a finitely generated
$R$-module. Suppose that $P$ is an ideal of $R$ such that $R/P$ is a
regular ring. Then $M$ has finite projective dimension if and only if
$\tor_i^R(M,R/P)=0$ for $i \gg 0$.
\end{lemma}

\begin{proof}
Since $R/P$ is regular, the residue field $K$ has a finite resolution
by free $R/P$-modules. Then $\tor_i^R(M,R/P)=0$ for $i \gg 0$ if and
only if $\tor_i^R(M,K)=0$ for $i \gg 0$. However $\tor_i^R(M,K)=0$ for 
$i \gg 0$ if and only if $M$ has finite projective dimension. 
\qed
\end{proof}

\section{Overview of the construction}
\label{overview}

We summarize the work that will be carried out in 
Sections~\ref{hypersurface module} and \ref{Gorenstein} and explain how this
provides the example we are aiming for.

In Section~\ref{hypersurface module} we construct a module $N$ of length $55$
and finite projective dimension over the local hypersurface
$$ 
A=K[U,V,W,X,Y,Z]_m/(UX+VY+WZ),
$$ 
where $m$ is the maximal ideal $(U,V,W,X,Y,Z)$, such that $N$ has a
nonzero intersection multiplicity with $A/P$ where $P$ is the prime
ideal $(u,v,w)$. Specifically, we have
$$
\chi(N,A/P) = \sum_{i=0}^5 (-1)^i \ell(\tor_i^A(N,A/P)) = -2.
$$ 
In Section~\ref{Gorenstein} we construct a Gorenstein normal
domain $R$ which is a module finite extension of $A$ and for which
there is an exact sequence of $A$-modules
$$
\begin{CD}
0 @>>> A^3 @>>> R @>>> P @>>> 0. {\textrm \qquad (**)}
\end{CD}
$$ 
Note that the ring $R$ has rank $4$ as an $A$-module. Consider the $R$-module 
$M=N\otimes_A R$. We claim that for this module
$$
\lim_{n\to\infty}\frac{\ell(F^n_R(M))}{p^{5n}} \neq \ell(M).
$$
To see this, let $\Fdot$ be a finite free resolution of $N$ over $A$. Since
$A \hookrightarrow R$ is a module-finite extension, the complex 
$\Gdot = \Fdot \otimes_A R$ has finite length homology. Furthermore since $R$ 
is Cohen-Macaulay, the complex $\Gdot$ is acyclic by the Acyclicity Lemma of 
Peskine and Szpiro, \cite[Lemma 1.8]{PS1}. Hence $\Gdot$ provides a finite free 
resolution of $M$ as an $R$-module. To compute the length of $M$ we use the 
additivity of $\chi(\Fdot\otimes_A -)$ on the exact sequence $(**)$. This gives
$$
\ell(M) = \chi(\Gdot) = \chi(\Fdot\otimes_A R) 
= 3\chi(\Fdot) + \chi(\Fdot\otimes_A P).
$$
The additivity also gives 
$$
\chi(\Fdot\otimes_A P) = \chi(\Fdot) - \chi(\Fdot\otimes_A A/P),
$$
and so
$$
\begin{array}{ll}
\ell(M) & = 4\chi(\Fdot) - \chi(\Fdot\otimes_A A/P) 
         = 4\ell(N) - \chi(N,A/P) \\ 
        & = 4\cdot 55 - (-2) = 222.
\end{array}
$$
On the other hand, since $R$ has rank $4$ as an $A$-module, Lemma~\ref{rank} 
gives
$$
\lim_{n\to\infty}\frac{\ell(F^n_R(M))}{p^{5n}} = 
\chinf(\Gdot) = \chinf(\Fdot\otimes_A R) 
= 4\chinf(\Fdot).
$$
Since $A$ is a hypersurface and $\Fdot$ is a finite resolution of $N$, we have 
$\chinf(\Fdot) = \ell(N)$ by \cite[Theorem 1.9]{Frobmult}. Therefore,
$$
\lim_{n\to\infty}\frac{\ell(F^n_R(M))}{p^{5n}} = 4\cdot 55 = 220.
$$ 

\section{A module of finite projective dimension}
\label{hypersurface module}

Consider the local hypersurface
$$
A=K[U,V,W,X,Y,Z]_m/(UX+VY+WZ)
$$ 
where $U,V,W,X,Y$ and $Z$ are indeterminates over a field $K$ of arbitrary 
characteristic, and $m$ is the maximal ideal $(U,V,W,X,Y,Z)$. We construct a 
module $N$ of
finite length and finite projective dimension over $A$, which has a nonzero
intersection multiplicity with the module $A/P$, where $P$ denotes the prime 
ideal $(u,v,w)A$. Specifically, we construct $N$ such that 
$$
\chi(N,A/P) = \sum_{i=0}^{5} (-1)^i \ell ({\tor}_i^A (N,A/P)) = -2.
$$

The following complex is an minimal free resolution of $A/P$:
$$
\begin{CD}
\cdots @>\phi_3>>A^4 @>\phi_4>>A^4 @>\phi_3>>A^4 @>\phi_2>>A^3 @>\phi_1>>A@>>>0.
\end{CD}
$$
The maps in this complex are given by the matrices:
$$
\phi_1 = \begin{pmatrix}
u & & v & & w \\
\end{pmatrix}, \ \ \ 
\phi_2 = \begin{pmatrix}
x &  &0   &-w  &  &v  \\
y &  &w   &0   &  &-u \\
z &  &-v  &u   &  &0  \\
\end{pmatrix},
$$
$$
\phi_3 = \begin{pmatrix}
0 &  & u  &  & v  &  & w  \\
u &  & 0  &  & z  &  & -y \\
v &  & -z &  & 0  &  & x  \\
w &  & y  &  & -x &  & 0  \\
\end{pmatrix}, \  \text{ and } \
\phi_4 = \begin{pmatrix}
0 &  & x  &  & y  &  & z  \\
x &  & 0  &  & -w &  & v  \\
y &  & w  &  & 0  &  & -u \\
z &  & -v &  & u  &  & 0  \\
\end{pmatrix}.
$$

The modules $\tor_i^A(N,A/P)$ may be computed by tensoring the above complex
with the module $N$. If $N$ has length $55$, the resulting complex may be
viewed as
$$
\begin{CD}
K^{220} @>\beta>>K^{220} @>\alpha>> K^{220} @>\theta_2>>K^{165} 
@>\theta_1>> K^{55}@>>> 0
\end{CD}
$$
where $\theta_1$, $\theta_2$, $\alpha$ and $\beta$ are matrices over
$K$. The $i$\,th homology of this complex is $\tor_i^A(N,A/P)$, and
will vanish for $i \ge 3$ provided the sum of the ranks of the
matrices $\alpha$ and $\beta$ is $220$. In this case the module $N$
has finite projective dimension by Lemma~\ref{projdim}, and an easy
calculation shows that
$$
\begin{aligned}
\chi(N &, A/P) \\
&= \ell(\tor_0^A(N,A/P)) - \ell(\tor_1^A(N,A/P)) + \ell(\tor_2^A(N,A/P)) \\
&= 55 - 165 + 220 - \rank(\alpha) = 110 - \rank(\alpha).
\end{aligned}
$$
In our construction, the matrix $\alpha$ will have rank $112$ and $\beta$
will have rank $108$.

As in \cite{DHM}, we regard a module of finite length over $A$ as a finite
dimensional vector space over $K$. The action of the generators of the ring can
then be treated as commuting nilpotent endomorphism of this vector space.  We
shall denote the endomorphisms given by the action of $u, v, w, x, y, z$  by
the matrices $\psi_1, \psi_2, \psi_3, \psi_4, \psi_5, \psi_6$, respectively.
Note that the matrices must satisfy the relation
$$
\psi_i\psi_j = \psi_j \psi_i, \text{ for all } i \text{ and }  j, 
$$
corresponding to commutativity, and the relation
$$ 
\psi_1 \psi_4 + \psi_2 \psi_5 + \psi_3 \psi_6 = 0,
$$
corresponding to the defining equation of the hypersurface.

The module of finite length and finite projection that we construct is
annihilated by $m^3 + (x,y,z)m$. Consequently $\psi_i$ may be written in block
form as
$$
\psi_i = \begin{pmatrix}
0 & & 0 & & a_i & & c_i \\
0 & & 0 & & 0   & & d_i \\
0 & & 0 & & 0   & & b_i \\
0 & & 0 & & 0   & & 0   \\
\end{pmatrix} \text{ for } i =1, \dots, 6.
$$
Furthermore we set 
$$
a_4 = a_5 = a_6 = 0 \text{ and } \ \ b_4 = b_5 = b_6 = 0.
$$
Since
$$
\psi_i \psi_j= \begin{pmatrix}
0 & & 0 & & 0  & & a_ib_j \\
0 & & 0 & & 0  & & 0      \\
0 & & 0 & & 0  & & 0      \\
0 & & 0 & & 0  & & 0      \\
\end{pmatrix}, 
$$
the relation $\psi_1 \psi_4 + \psi_2 \psi_5 + \psi_3 \psi_6 = 0$ is easily seen 
to be satisfied. 

The matrices $a_1, a_2, a_3$ are
$$
a_1 = \begin{pmatrix}
1  & & 0 & & 0 \\
\end{pmatrix}, \ \
a_2 = \begin{pmatrix}
0  & & 1 & & 0 \\
\end{pmatrix}, \ \
a_3 = \begin{pmatrix}
0  & & 0 & & 1 \\
\end{pmatrix},
$$
where $1$ denotes the $4 \times 4$ identity matrix, and $0$ denotes the
$4 \times 4$ zero matrix. For $b_i$ we take
$$
b_1 = \begin{pmatrix}
1  &0 &0 &0 &0 &0 \\
0  &1 &0 &0 &0 &0 \\
0  &0 &1 &0 &0 &0 \\
\end{pmatrix}, \ \
b_2 = \begin{pmatrix}
0  &1 &0 &0 &0 &0 \\
0  &0 &0 &1 &0 &0 \\
0  &0 &0 &0 &1 &0 \\
\end{pmatrix}, \ \
b_3 = \begin{pmatrix}
0  &0 &1 &0 &0 &0 \\
0  &0 &0 &0 &1 &0 \\
0  &0 &0 &0 &0 &1 \\
\end{pmatrix}
$$
where $1$ denotes the $4 \times 4$ identity matrix, and $0$ denotes the
$4 \times 4$ zero matrix. Note that for all $i$ and $j$ we have 
$a_ib_j = a_jb_i$, and so the commutativity relations 
$\psi_i\psi_j = \psi_j \psi_i$ do hold.

After interchanging certain rows and columns, the matrix $\alpha$ is
$$
\begin{pmatrix}
0    & a_1  & a_2  & a_3  & 0   & c_1  & c_2  & c_3  \\
a_1  & 0    & 0    & 0    & c_1 & 0    & c_6  & -c_5 \\
a_2  & 0    & 0    & 0    & c_2 & -c_6 & 0    & c_4  \\
a_3  & 0    & 0    & 0    & c_3 & c_5  & -c_4 & 0    \\
0    & 0    & 0    & 0    & 0   & d_1  & d_2  & d_3  \\
0    & 0    & 0    & 0    & d_1 & 0    & d_6  & -d_5 \\
0    & 0    & 0    & 0    & d_2 & -d_6 & 0    & d_4  \\
0    & 0    & 0    & 0    & d_3 & d_5  & -d_4 & 0    \\
0    & 0    & 0    & 0    & 0   & b_1  & b_2  & b_3  \\
0    & 0    & 0    & 0    & b_1 & 0    & 0    & 0    \\
0    & 0    & 0    & 0    & b_2 & 0    & 0    & 0    \\
0    & 0    & 0    & 0    & b_3 & 0    & 0    & 0    \\
\end{pmatrix}.
$$
The rank of the matrix $\alpha$ is easily seen to be the sum of $40$
and the rank of the submatrix
$$ \alpha_1 = 
\begin{pmatrix}
d_1  & d_2  & d_3  \\
0    & d_6  & -d_5 \\
-d_6 & 0    & d_4  \\
d_5  & -d_4 & 0    \\
b_1  & b_2  & b_3  \\
\end{pmatrix}.
$$

Similarly, after deleting rows and columns of zeros and interchanging certain
rows and columns, the matrix $\beta$ reduces to
$$
\begin{pmatrix}
0    & 0    & 0    & 0   & c_4  & c_5  & c_6  \\
0    & -a_3 & a_2  & c_4 & 0    & -c_3 & c_2  \\
a_3  & 0    & -a_1 & c_5 & c_3  & 0    & -c_1 \\
-a_2 & a_1  & 0    & c_6 & -c_2 & c_1  & 0    \\
0    & 0    & 0    & 0   & d_4  & d_5  & d_6  \\
0    & 0    & 0    & d_4 & 0    & -d_3 & d_2  \\
0    & 0    & 0    & d_5 & d_3  & 0    & -d_1 \\
0    & 0    & 0    & d_6 & -d_2 & d_1  & 0    \\
0    & 0    & 0    & 0   & 0    & -b_3 & b_2  \\
0    & 0    & 0    & 0   & b_3  & 0    & -b_1 \\
0    & 0    & 0    & 0   & -b_2 & b_1  & 0    \\
\end{pmatrix}.
$$
This rank of this matrix is the sum of $12$ and the rank of the submatrix
$$ \beta_1 =
\begin{pmatrix}
0   & c_4  & c_5  & c_6  \\
0   & d_4  & d_5  & d_6  \\
d_4 & 0    & -d_3 & d_2  \\
d_5 & d_3  & 0    & -d_1 \\
d_6 & -d_2 & d_1  & 0    \\
0   & 0    & -b_3 & b_2  \\
0   & b_3  & 0    & -b_1 \\
0   & -b_2 & b_1  & 0    \\
\end{pmatrix}.
$$

We next let $c_1=c_2=c_3=c_4=c_5=0$ and 
$$
c_6=\begin{pmatrix}
0 & 0 & 0 & 0 & 0 & 1 \\
\end{pmatrix}
$$
where $1$ denotes the $4 \times 4$ identity matrix, and $0$ denotes the
$4 \times 4$ zero matrix. It remains to exhibit matrices $d_i$ for $1 \le i \le
6$ such that the matrices 
$\alpha$ and $\beta$ have ranks $112$ and $108$ respectively. We let 
$d_1=d_2=0$ and 
$$d_3=\left(
\begin{array}{llllllllllllllllllllllll}
0&0&0&0&0&0&0&0&1&0&0&0&0&0&0&0&0&0&0&0&0&0&0&0 \\
0&0&0&0&0&0&0&0&0&1&0&0&0&0&0&0&0&0&0&0&0&0&0&0 \\
0&0&0&0&0&0&0&0&0&0&1&0&0&0&0&0&0&0&0&0&0&0&0&0 \\
0&0&0&0&0&0&0&0&0&0&0&1&0&0&0&0&0&0&0&0&0&0&0&0 \\
0&0&0&0&0&0&0&0&0&0&0&0&0&0&0&0&1&0&0&0&0&0&0&0 \\
0&0&0&0&0&0&0&0&0&0&0&0&0&0&0&0&0&1&0&0&0&0&0&0 \\
0&0&0&0&0&0&0&0&0&0&0&0&0&0&0&0&0&0&1&0&0&0&0&0 \\
0&0&0&0&0&0&0&0&0&0&0&0&0&0&0&0&0&0&0&1&0&0&0&0 \\
0&0&0&0&0&0&0&0&0&0&0&0&0&0&0&0&0&0&0&0&1&0&0&0 \\
0&0&0&0&0&0&0&0&0&0&0&0&0&0&0&0&0&0&0&0&0&1&0&0 \\
0&0&0&0&0&0&0&0&0&0&0&0&0&0&0&0&0&0&0&0&0&0&1&0 \\
0&0&0&0&0&0&0&0&0&0&0&0&0&0&0&0&0&0&0&0&0&0&0&1 \\
0&0&0&0&0&0&0&0&0&0&0&0&1&0&0&0&0&0&0&0&0&0&0&0 \\
0&0&0&0&0&0&0&0&0&0&0&0&0&1&0&0&0&0&0&0&0&0&0&0 \\
0&0&0&0&0&0&0&0&0&0&0&0&0&0&1&0&0&0&0&0&0&0&0&0 \\
\end{array}
\right),
$$

$$d_4=\left(
\begin{array}{llllllllllllllllllllllll}
1&0&0&0&0&0&0&0&0&0&0&0&0&0&0&0&0&0&0&0&1&0&0&0 \\
0&1&0&0&0&0&0&0&0&0&0&0&0&0&0&0&0&0&0&0&0&1&0&0 \\
0&0&1&0&0&0&0&0&0&0&0&0&0&0&0&0&0&0&0&0&0&0&1&0 \\
0&0&0&1&0&0&0&0&0&0&0&0&0&0&0&0&0&0&0&0&0&0&0&1 \\
0&0&0&0&1&0&0&0&0&0&0&0&0&0&0&0&0&0&0&0&0&0&0&0 \\
0&0&0&0&0&1&0&0&0&0&0&0&0&0&0&0&0&0&0&0&0&0&0&0 \\
0&0&0&0&0&0&1&0&0&0&0&0&0&0&0&0&0&0&0&0&0&0&0&0 \\
0&0&0&0&0&0&0&1&0&0&0&0&0&0&0&0&0&0&0&0&0&0&0&0 \\
0&0&0&0&0&0&0&0&0&0&0&0&0&0&0&0&0&0&0&0&0&0&0&0 \\
0&0&0&0&0&0&0&0&0&0&0&0&0&0&0&0&0&0&0&0&0&0&0&0 \\
0&0&0&0&0&0&0&0&0&0&0&0&0&0&0&0&0&0&0&0&0&0&0&0 \\
0&0&0&0&0&0&0&0&0&0&0&0&0&0&0&0&0&0&0&0&0&0&0&0 \\
0&0&0&0&0&0&0&0&0&0&0&0&0&0&0&0&0&0&0&0&0&0&0&0 \\
0&0&0&0&0&0&0&0&0&0&0&0&1&0&0&0&0&0&0&0&0&0&0&0 \\
0&0&0&0&0&0&0&0&0&0&0&0&0&0&0&1&0&0&0&0&0&0&0&0 \\
\end{array}
\right),
$$

$$d_5=\left(
\begin{array}{llllllllllllllllllllllll}
0&0&0&0&0&0&0&0&0&0&0&0&0&0&0&0&0&0&0&0&0&0&0&0 \\
0&0&0&0&0&0&0&0&0&0&0&0&0&0&0&0&0&0&0&0&0&0&0&0 \\
0&0&0&0&0&0&0&0&0&0&0&0&0&0&0&0&0&0&0&0&0&0&0&0 \\
0&0&0&0&0&0&0&0&0&0&0&0&0&0&0&0&0&0&0&0&0&0&0&0 \\
0&0&0&0&0&0&0&0&0&0&0&0&0&0&0&0&1&0&0&0&0&0&0&0 \\
0&0&0&0&0&0&0&0&0&0&0&0&0&0&0&0&0&1&0&0&0&0&0&0 \\
0&0&0&0&0&0&0&0&0&0&0&0&0&0&0&0&0&0&1&0&0&0&0&0 \\
0&0&0&0&0&0&0&0&0&0&0&0&0&0&0&0&0&0&0&1&0&0&0&0 \\
0&0&0&0&0&0&0&0&0&0&0&0&0&0&0&0&0&0&0&0&1&0&0&0 \\
0&0&0&0&0&0&0&0&0&0&0&0&0&0&0&0&0&0&0&0&0&1&0&0 \\
0&0&0&0&0&0&0&0&0&0&0&0&0&0&0&0&0&0&0&0&0&0&1&0 \\
0&0&0&0&0&0&0&0&0&0&0&0&0&0&0&0&0&0&0&0&0&0&0&1 \\
0&0&0&0&0&0&0&0&0&0&0&0&0&0&0&1&0&0&0&0&0&0&0&0 \\
0&0&0&0&0&0&0&0&0&0&0&0&0&0&0&0&0&0&0&0&0&0&0&0 \\
0&0&0&0&0&0&0&0&0&0&0&0&0&0&0&0&0&0&0&0&0&0&0&0 \\
\end{array}
\right),
$$
and

$$d_6=\left(
\begin{array}{llllllllllllllllllllllll}
1&0&0&0&0&0&0&0&0&0&0&0&0&0&0&0&0&0&0&0&0&0&0&0 \\
0&1&0&0&0&0&0&0&0&0&0&0&0&0&0&0&0&0&0&0&0&0&0&0 \\
0&0&1&0&0&0&0&0&0&0&0&0&0&0&0&0&0&0&0&0&0&0&0&0 \\
0&0&0&1&0&0&0&0&0&0&0&0&0&0&0&0&0&0&0&0&0&0&0&0 \\
0&0&0&0&1&0&0&0&0&0&0&0&0&0&0&0&0&0&0&0&0&0&0&0 \\
0&0&0&0&0&1&0&0&0&0&0&0&0&0&0&0&0&0&0&0&0&0&0&0 \\
0&0&0&0&0&0&1&0&0&0&0&0&0&0&0&0&0&0&0&0&0&0&0&0 \\
0&0&0&0&0&0&0&1&0&0&0&0&0&0&0&0&0&0&0&0&0&0&0&0 \\
0&0&0&0&0&0&0&0&1&0&0&0&0&0&0&0&0&0&0&0&0&0&0&0 \\
0&0&0&0&0&0&0&0&0&1&0&0&0&0&0&0&0&0&0&0&0&0&0&0 \\
0&0&0&0&0&0&0&0&0&0&1&0&0&0&0&0&0&0&0&0&0&0&0&0 \\
0&0&0&0&0&0&0&0&0&0&0&1&0&0&0&0&0&0&0&0&0&0&0&0 \\
0&0&0&0&0&0&0&0&0&0&0&0&1&0&0&0&0&0&0&0&0&0&0&0 \\
0&0&0&0&0&0&0&0&0&0&0&0&0&1&0&0&0&0&0&0&0&0&0&0 \\
0&0&0&0&0&0&0&0&0&0&0&0&0&0&1&0&0&0&0&0&0&0&0&0 \\
\end{array}
\right).
$$
Of course, any choice of matrices $d_3, d_4, d_5$ and $d_6$ in general position
will serve our purpose. That the matrices exhibited above do have the
required rank properties is an elementary, though tedious, verification. 

\section{Construction of Gorenstein rings}
\label{Gorenstein}

We now construct a Gorenstein normal domain $R$ which is an extension of a
free $A$-module by the prime ideal $P$. The construction is carried out in the
case that $A$ is a hypersurface over a field of characteristic $2$.

Consider the ring $R = A[a,b,c,d,e]$ where 
$$
a=\sqrt{uyz}, \ b=\sqrt{vxz}, \ c=\sqrt{wxy}, \ d=\sqrt{uvw}, \ e=\sqrt{vwyz}.
$$
The ring $R$ is a normal domain; in fact, it is the integral closure
of $A$ in the field $L( \sqrt{uyz}, \ \sqrt{vxz})$ where $L$ is the
fraction field of $A$. This furthermore shows that the ring $R$ has 
rank $4$ as an $A$-module. 

To show that $R$ is a Cohen-Macaulay ring, we work with the system of
parameters $u,x,v,y,w-z$ and show that the multiplicity of the ideal 
$I=(u,x,v,y,w-z)$ is $8$ and that this equals the length of the
$K$ vector space $R/I$.

Consider the extension $K[u,x,v,y,w-z] \subseteq R$. The degree of this
extension may be computed by examining the corresponding extension of
fraction fields, and so is easily seen to be $8$. 

The following relations show that $R$ is generated as an $A$-module
by the elements $1, a, b, c, d, e$. 

$$
\begin{array}{lllll}
a^2 = uyz, & &
ab  = ze+vyz, & &  
ac  = ye+wyz, \\
ad  = ue, & &
ae  = dyz, & &
b^2  = vxz, \\
bc  = xe, & &
bd  = ve+vwz, & & 
be  = vzc, \\
c^2 = wxy, & &
cd  = we+vwy, & & 
ce  = wyb, \\ 
d^2 = uvw, & &
de  = vwa, & &
e^2 = vwyz.
\end{array}
$$
The relations amongst the elements of $R$ also include
$$
\begin{array}{llllll}
ux+vy+wz, & &
xa+yb+zc, & &
wb+vc+xd, \\
wa+uc+yd, & &
va+ub+zd. & &
\end{array}
$$
Consequently the images of the following elements form a generating 
set for $R/I$ as a $K$ vector space:
$$
1, \ z, \ a, \ b, \ c, \ d, \ e, \ ze.
$$
Hence the length of $R/I$ is less than or equal to $8$, but since the
multiplicity of the ideal $I$ was earlier computed to be $8$, the
length must be precisely $8$. This shows that $R$ is a Cohen-Macaulay
ring, and furthermore that it is Gorenstein, since the socle in $R/I$
is of dimension one, being generated by the image of $ze$. It is not
difficult to verify that the relations listed above are precisely the
relations amongst the generators of $R$.

Lastly, it may be verified that there is an exact sequence of $A$-modules
$$
\begin{CD}
0 @>>> A^3 @>>> R @>>> P @>>> 0
\end{CD}
$$
where the map $R \to P$ is determined by 
$$
1 \to 0, \ a \to u, \ b \to v, \ c \to w, \ d \to 0, \ e \to 0.
$$

\section{Local Chern Characters}
\label{Chern}

Let $X$ be a closed subset of $\Spec R$ where the ring $R$ is a
$d$-dimensional homomorphic image of a regular local ring. For each
integer $i$, let $Z_i(X)$ denote the free $\mathbb Q$-module generated
by cycles of the form $[R/P]$ for $P$ in $X$ such that $\dim R/P = i$.
For a prime $Q$ with $\dim R/Q = i+1$ and an element $x$ of $R$ such
that $x \notin Q$, define
$$
\div(x,Q) = \sum_{\dim R/P = i} \ell((R/(Q+xR))_P) [R/P].
$$
Let $B_i(X)$ denote the subgroup of $Z_i(X)$ generated by elements of
the form $\div(x,Q)$ for $Q$ in $X$. The $i$\,th graded piece of the
Chow group of $X$ is $A_i(X) = Z_i(X)/B_i(X)$, and the Chow group of
$X$ is
$$
A_*(X) = \bigoplus_{i=1}^d A_i(X). 
$$

Let $\Gdot$ be a bounded free complex of $R$-modules, and let
$Z\subset \Spec R$ be the support of this complex. The local Chern
character $\ch(\Gdot)$ is the sum
$$
\ch(\Gdot) = \ch_d(\Gdot) + \ch_{d-1}(\Gdot) + \cdots + \ch_0(\Gdot)
$$ 
where, for each $k$, 
$$
\ch_k(\Gdot) : A_i(X) \to A_{i-k}(X \cap Z)
$$
is a $\mathbb Q$-module homomorphism. (For precise definitions and
properties we refer the reader to \cite{Fu}.) We consider the special
case of the local Riemann-Roch formula where the homology modules of
the complex $\Gdot$ are of finite length: for any finitely generated
$R$-module $N$, there is an element
$$
\tau(N) = \tau_d(N) + \cdots + \tau_0(N) \in A_*(\Supp(N)),
$$ 
called the {\it Todd class}\/ of $N$, such that 
$$
\chi(\Gdot\otimes N) = \ch(\Gdot)(\tau(N)).
$$
Note that since the support of $\Gdot$ is the closed point of $\Spec R$, 
$\ch(\Gdot)(\tau(N))$ is an element of $A_{0}(\Spec (R/m)) \cong \mathbb Q.$
If $\Gdot$ is the resolution of a module $M$ of finite length and finite
projective dimension, the local Riemann-Roch formula then gives  
$$
\ell(M) = \ch(\Gdot)(\tau(R)) = \sum_{i=0}^d \ch_i(\Gdot)(\tau_i(R)).
$$

Now suppose in addition that $R$ is a complete local  ring over a perfect
field of prime characteristic $p$. Since the local Chern characters are
compatible with finite maps, one can show that
$$
\lim_{n\to\infty}\frac{\ell(F^n_R(M))}{p^{nd}} = \ch_d(\Gdot)(\tau_d(R)),
$$ 
see, \cite[Proposition 12.7.1]{RobBook}. If the ring $R$
is a complete intersection then $\tau_i(R)= 0$ for all $i < d$, and so
$$
\ell(M) = \ch(\Gdot)(\tau(R)) = \ch_d(\Gdot)(\tau_d(R)) 
=\lim_{n\to\infty}\frac{\ell(F^n_R(M))}{p^{nd}}.
$$ 
For a Gorenstein ring $R$ the duality property shows that $\tau_{d-i}(R) = 0$ 
for all odd numbers $i$. In the specific case that $R$ is a Gorenstein ring of 
dimension $3$ and $\Gdot$ is the resolution of a module $M$ of finite length and 
finite projective dimension, this then gives
$$
\ell(M) = \ch(\Gdot)(\tau(R)) =\ch_3(\Gdot)(\tau_3(R)) +\ch_1(\Gdot)(\tau_1(R)).
$$
The operator $\ch_1(\Gdot)$ can be identified with the MacRae
invariant, and is known to vanish whenever the module $M$ has
codimension greater than one, see \cite[Theorem 3]{RobMacRae}, and
also \cite{Foxby,MR}. In \cite{Kurano} Kurano gave examples of
Gorenstein rings $R$ for which $\tau_i(R)\neq 0$ for some
$i<d$. However it is not known in the case of Kurano's examples if
there exists a bounded free complex $\Gdot$ with finite length
homology for which $\ch_i(\Gdot)(\tau_i(R))$ does not vanish for some
$i<d$.  We would also like to point out that, using the example in
\cite{DHM} of a module with negative intersection multiplicity,
Roberts did construct an example of a \CM ring $R$ of dimension $3$
such that $\ch_2(\Gdot)(\tau_2(R))\neq 0$ for a bounded free complex
$\Gdot$ with finite length homology, which is, in fact, the resolution
of a module of finite length and finite projective dimension.

In our example, where $R$ is a Gorenstein ring of dimension $5$ and $\Gdot$ is
the finite free resolution of the module $M$, we have 
$$
\ell(M) = \ch_5(\Gdot)(\tau_5(R)) + \ch_3(\Gdot)(\tau_3(R)) = 222.
$$
(Recall that $\ch_1(\Gdot)=0$ by \cite[Theorem 3]{RobMacRae}.) On the
other hand,
$$
\ch_5(\Gdot)(\tau_5(R)) = \lim_{n\to\infty}\frac{\ell(F^n_R(M))}{p^{5n}} = 220,
$$ 
and so, we must have
$$
\ch_3(\Gdot)(\tau_3(R)) = 2.
$$ 
Thus our example provides an example of a Todd class $\tau_3(R)$ over
a five dimensional Gorenstein ring $R$, and of a bounded free complex
with a local Chern character that does not vanish on this class.
Furthermore, we may also conclude that in our example, where $p=2$, we
have the formula
$$
\ell(F^n_R(M)) = 220\cdot2^{5n} + 2\cdot2^{3n} \text{ \ for all } n \ge 0
$$
by the following proposition, which follows from the compatibility of
local Chern characters with finite maps: see the proof of 
\cite[Proposition 12.7.1]{RobBook}.

\begin{proposition}[Roberts]
Let $R$ be a Noetherian local ring of positive characteristic $p$ and
dimension $d$. Suppose that the residue field of $R$ is perfect and
that the Frobenius endomorphism is a finite map. Then for any bounded
free complex $\Gdot$ with finite length homology, we have
$$
\ch(F^n_R(\Gdot))(\tau(R)) = \sum_{i=0}^d p^{ni}\ch_i(\Gdot)(\tau_i(R))
$$
for all $n\geq 0$.
\end{proposition}

\section{Further consequences of the example}
\label{further}

Let $(R,m)$ be a complete local ring of dimension $d$, and let $M$ be an
$R$-module of finite length. The length $\ell(F^n(M))$ is equal to 
$p^{nd}\ell(M)$ if the ring $R$ is regular, and this may be viewed as a special
case of the following result, \cite[Theorem 1.9]{Frobmult}:

\begin{theorem}[Dutta]
Let $(R,m)$ be a complete intersection ring of dimension $d$, and let $M$ be an
$R$-module of finite length. Then 
$$
\ell(F^n_R(M)) \geq \ell(M) p^{nd},
$$
and equality holds if the module $M$ has finite projective dimension.
\end{theorem}

C.~Miller has obtained a converse to this theorem in \cite{CM}, which then
gives a characterisation of modules of finite length and finite projective
dimension:

\begin{theorem}[Miller]
Let $(R,m)$ be a complete intersection ring of dimension $d$, and let $M$ be an
$R$-module of finite length. Then the following statements are equivalent:
\item $1)$ \quad $\ell(F^n_R(M)) = \ell(M) p^{nd}$ for all $n \ge 0$.
\item $2)$ \quad The module $M$ has finite projective dimension.
\item $3)$ \quad $\lim_{n \to \infty}\frac{\ell(F^n_R(M))}{p^{nd}} = \ell(M)$.
\end{theorem}

A natural question raised by these theorems is whether there is a relationship
between $\ell(F^n_R(M))$ and $\ell(M) p^{nd}$ if the ring $R$ is Gorenstein, but
is not a complete intersection. The examples we have constructed already show
that over a Gorenstein ring $R$ of dimension $5$, the equality $\ell(F^n_R(M)) =
\ell(M) p^{5n}$ need not hold for a module $M$ of finite length and finite
projective dimension. We next show that the inequality $\ell(F^n_R(M)) \geq
\ell(M) p^{nd}$ which holds for {\it all}\/ modules $M$ of finite length over a
complete intersection ring, fails to hold over  Gorenstein rings even when the
finite length module $M$ has finite projective dimension.

In Section~\ref{hypersurface module} we constructed a module $N$ over the 
hypersurface
$$ 
A=K[U,V,W,X,Y,Z]_m/(UX+VY+WZ) 
$$ 
which has finite length and finite projective dimension and has the property 
that 
$$
\chi(N,A/P) =  -2,
$$
where $P$ is the prime ideal $(u,v,w)$.  Let $Q$ denote the prime ideal
$(x,v,w)$, and consider the short exact sequence
$$
\begin{CD}
0 @>>> A/Q @>>> A/(v,w) @>>> A/P @>>> 0,
\end{CD}
$$
where the first map is multiplication by the element $u$.  Since the elements 
$v$ and $w$ 
form a regular sequence in the ring $A$, the Koszul resolution gives that 
$$
\chi(N,A/(v,w)) = \sum_{i=0}^2 (-1)^i H_i(v,w;N),
$$
where $H_i(v,w;N)$ denotes the $i \,$th Koszul homology module.
Since $\dim (N) < 2$, we then have $\sum_{i=0}^2 (-1)^i H_i(v,w;N) =0$
by \cite[Theorem 1, Chapter IV]{Serre}. Consequently, 
$$
\chi(N,A/Q) = -\chi(N,A/P) = 2.
$$
Consider the automorphism $\sigma$ of $A$ which switches $x$ and $u$, and fixes 
$v$, $w$, $y$ and $z$. Let $N'$ denote the module $N$ now viewed as an 
$A$-module by restriction of scalars via $\sigma$. We then have
$$
\chi(N',A/P) = \chi(N,A/Q) = 2.
$$
The same argument as in Section~\ref{overview} then shows that the
$R$-module $M' = N' \otimes_A R$ has length
$$
\ell(M') = 4\ell(N') - \chi(N',A/P)  = 218,
$$
whereas
$$
\lim_{n\to\infty}\frac{\ell(F^n_R(M'))}{p^{5n}} = 220.
$$ 

\begin{remark}
For a \CM local ring $A$, let \AE(A) denote the Grothen\-dieck group of
$A$-modules which have finite length and finite projective dimension. If 
$x_1, \dots,
x_d$ is a system of parameters for $A$, then $[A/(x_1, \dots, x_d)]$ is an
element of \AE(A). The intersection theory of the modules of the form $A/(x_1,
\dots, x_d)$ is well understood: if $M$ is a finitely generated $A$-module
with $\dim M < d$, then  
$$ 
\chi(A/(x_1, \dots, x_d), M) = 0 
$$ 
by \cite[Theorem 1, Chapter IV]{Serre}. Consequently it is of interest
to understand the group $G(A)$ which is the quotient of \AE(A) by the
subgroup generated by all modules of the form $A/(x_1, \dots, x_d)$
where $x_1, \dots, x_d$ is a system of parameters for $A$. Using
$K$-theoretic methods, M.~Levine has computed this group for certain
hypersurfaces, see \cite[\S 4]{Levine}. We would like to point out
that in the case where $A$ is the hypersurface
$$
A=K[U,V,W,X,Y,Z]_m/(UX+VY+WZ),
$$
and $N$ is the module of finite length and finite projective dimension
constructed in Section~\ref{hypersurface module}, the fact that
$\chi(N,A/(u,v,w))\neq 0$ shows that $[N]$ is a nonzero element of the
group $G(A)$.
\end{remark}

\begin{acknowledgement}
We thank Mel Hochster for valuable discussions and encouragement
during the course of this work. 
\end{acknowledgement}

\end{document}